\documentclass[12pt,leqno]{amsart}
\usepackage[dvips]{graphicx}
\usepackage{amsfonts}
\usepackage{amssymb}
\usepackage{amsmath}
\usepackage{amscd}
\usepackage{graphicx}
\def\DATE{ 4 Chevat 5777 }
\newtheorem{theorem}{Theorem}
\newtheorem{definition}[theorem]{Definition}
\newtheorem{corollary}[theorem]{Corollary}

\newtheorem{proposition}[theorem]{Proposition}

\newcommand\C{\mathbb{C}}

\newcommand\K{\mathbb{K}}

\newcommand\ds{\displaystyle}
\newcommand\pf{\noindent{\it Proof. }}

\email{elisabeth.remm@uha.fr ; goze.rac@gmail.com}

\begin{document}
\title{On the algebraic variety of Hom-Lie algebras}
\author{Elisabeth Remm, Michel Goze}
  \date{ 4 Chevat 5777 }
 \maketitle

% \thanks{$^\ast$Corresponding author. Email: elisabeth.remm@uha.fr} \\

\begin{abstract}
The set $\mathcal{H}Lie_n$ of the $n$-dimensional Hom-Lie algebras over an algebraically closed field of characteristic zero is provided with a structure of algebraic subvariety of the affine plane $\K^{n^2(n-1)/2}$. For $n=3$, these two sets coincide, for $n=4$ it is an hypersurface in $\K^{24}$. For $n \geq 5$, we describe the scheme of polynomial equations which define $\mathcal{H}Lie_n$. We determine also what are the classes of Hom-Lie algebras which are $\mathcal{P}$-algebras where $\mathcal{P}$ is a binary quadratic operads.
  \end{abstract}


\medskip


\medskip

\section{Introduction}

The notion of Hom-Lie algebras was introduced by Hartwig, Larsson and  Silvestrov  in \cite{sil}. Their principal motivation concerns deformation of the Witt algebra. This Lie algebra is  the complexification of the Lie algebra of polynomial vector fields on a circle. A basis for the Witt algebra is given by the vector fields 
$$L_n=-z^{n+1}\ds \frac{\partial}{\partial z}$$
 for any  $n \in \mathbb {Z}.$
The Lie bracket  is given by
$$[L_{m},L_{n}]=(m-n)L_{{m+n}}.$$
The Witt algebra is also viewed as  the Lie algebra of derivations of the ring $\C[z,z^{-1}]$. Recall that a derivation on an algebra with product denoted by $ab$ is a linear operator satisfying $D(ab)=D(a)b+aD(b)$. The Lie bracket of two derivations $D$ and $D'$ is  $[D,D']=D \circ D' - D' \circ D.$ We can also define a new class of linear operators generalizing derivations, the Jackson derivate, given by
$$D_q(f)(z)=\ds \frac{f(qz)-f(z)}{qz-z}.$$
It is clear that $D_q$ is a linear operator, but its behavior on the product is quite different as the classical derivative:
$$D_q(fg(z))=g(z)D_q(f(z))+f(qz)D_q(g(z)).$$
The authors interpret this relation by putting
\begin{equation}
\label{sigma}
D_q(fg)=gD_q(f)+\sigma (f)D_q(g)
\end{equation}
where $\sigma$  is given by $\sigma(f)(z)=f(qz)$ for any $f \in \C[z,z^{-1}]$. Starting from (\ref{sigma}) and for a given $\sigma$, one defines a new space of  "derivations" on  $C[z,z^{-1}]$ constituted of linear operator $D$ satisfying this relation. With the classical bracket we obtain is a new type of algebra so called $\sigma$-deformations of the Witt algebra.  This new approach  leads  naturally to considerer the space of $\sigma$-derivations, that is, linear operators satisfying (\ref{sigma}), to provide it with the multiplication associated with a bracket. This new algebra is not a Lie algebra because the bracket doesn't satisfies the Jacobi conditions. The authors shows that this bracket satisfies a "generalized Jacobi condition". They have called this new class of algebras the class of Hom-Lie algebras. This notion, introduced in \cite{sil},
 since made the object of numerous studies and was also generalized (see \cite{MDL}).
We denote by $\mathcal{SA}lg_n$ the set of the $n$-dimensional skew-symmetric algebras  $(\mathcal{A},\mu)$ over an algebraically closed field $\K$ of characteristic $0$ whose multiplication $\mu$ is skew-symmetric and by $\mathcal{H}Lie_n$ the subset of Hom-Lie algebras. It is clear that $\mathcal{SA}lg_n$ is an affine variety isomorphic to $K^{n^2(n-1)/2}$. In this work, we show that $\mathcal{H}Lie_n$ is an affine algebraic subvariety of $\mathcal{SA}lg_n$. We study in particular the case $n=3$ and $n=4$, proving that in dimension $3$ any skew-symmetric algebra is a Hom-Lie algebra and in dimension $4$,  $\mathcal{H}Lie_4$ is an algebraic hypersurface in $\mathcal{SA}lg_4$. We end this work by the determination of binary quadratic operads whose associated algebras are Hom-Lie.

\section{Hom-Lie algebras \cite{sil}}

Before to begin the study of Hom-Lie algebras, we recall some very classical results on the notions of algebras and fix the used notations.

An algebra  over a field $\K$ of characteristic $0$ is a pair $(\mathcal{A},\mu)$ where $\mathcal{A}$ is $\K$-vector space and $\mu$ a bilinear map
$$\mu : \mathcal{A} \times \mathcal{A} \rightarrow \mathcal{A}$$
which is classically called the multiplication of $\mathcal{A}$. When any confusion is possible, we shall denote it $\mathcal{A}$ in place of $(\mathcal{A},\mu)$. Two algebras $(\mathcal{A},\mu)$ and $(\mathcal{A}',\mu')$ are isomorphic (in terms of algebras), if there is a linear isomorphism
$$f :\mathcal{A} \rightarrow \mathcal{A}'$$
satisfying
\begin{equation}
\label{iso}
f(\mu(X,Y))=\mu'(f(X),f(Y))
\end{equation}
for any $X,Y \in \mathcal{A}$.

We assume now that all considered algebras are of finite dimension. Let $\{e_1,\cdots,e_n\}$ be a basis of the algebra $(\mathcal{A},\mu)$. The constants $C_{i,j}^k$ given by
$$\mu(e_i,e_j)=\sum_{k}C_{i,j}^k e_k$$
are called the structure constants of $\mu$ (or $\mathcal{A}$)  relative to the basis $\{e_1,\cdots,e_n\}$.  The set $\mathcal{A}lg_n$ of $n$-dimensional $\K$-algebras can be identified to the vector space $\K^{n^3}$, identifying the multiplication $\mu$ with its structure constants $C_{i,j}^k$ related to a given basis of $\K^n$. One of the main problem is to give the classification of the algebras of a given dimension. This problem was solved for the dimensions $2$ and $3$ (in this last case for the skew-symmetric algebras). This classification corresponds to the determination of all orbits related to the action (\ref{iso}) of the linear group $GL(n,\K)$ on $\mathcal{A}lg_n.$

Among all the algebras, certain classes play importing roles, as the class of  Lie algebras, associative algebras, pre-Lie algebras. The aim of this work concerns a recent class, the Hom-Lie algebras.

\begin{definition}
A Hom-Lie algebra is a triple $(A, \mu, \alpha)$
consisting of a linear space $A$, a skew-bilinear map $\mu: V\times V \rightarrow V$ and
 a linear space homomorphism $f: A \rightarrow A$
 satisfying the Hom-Jacobi identity
$$
 \circlearrowleft_{x,y,z}{\mu(\mu(x,y),f(z))}=0
$$
for all $x, y, z$ in $A$, where
$\circlearrowleft_{x,y,z}$ denotes summation over
the cyclic permutations on $x,y,z$.
\end{definition}

For example, a Hom-Lie algebra whose endomorphism $f$ is the identity is a Lie algebra. We deduce, 
since any $2$-dimensional skew-symmetric algebra (the multiplication $\mu$ is a skew-symmetric bilinear map)  is a Lie algebra, that any $2$-dimensional algebra is a Hom-Lie algebra.

In the following section we are interested by the determination of all Hom-Lie algebras for small dimensions.

\section{Hom-Lie algebras of dimension $3$}

We use the original approach developped in \cite{ERdim3}. Any  $3$-dimensional $\K$-algebra $(\mathcal{A},\mu)$ is defined by its  structure constants $\{\alpha_i,\beta_i, \gamma_i \}_{i=1,2,3}$ with respect to a  given basis $\{e_1,e_2,e_3\}:$ 
$$
\left\{
\begin{array}{l}
\mu(e_1,e_2)=\alpha_1e_1+\beta _1e_2+ \gamma_1e_3, \\
\mu(e_1,e_3)=\alpha_2e_1+\beta _2e_2+\gamma_2e_3, \\
\mu(e_2,e_3)=\alpha_3e_1+\beta _3e_2+\gamma_3e_3. \\
 \end{array}
 \right.
 $$
Let $f$ be an element of  $gl(3,\K)$ and  consider its matrix in the same basis $\{e_1,e_2,e_3\}$
$$\left(
\begin{array}{lll}
   a_1 & b_1 & c_1   \\
     a_2 & b_2 & c_2   \\
   a_3 & b_3 & c_3   \\
\end{array}
\right).
$$
We then define the vector
$$v_f=\ ^t (a_1,a_2,a_3,b_1,b_2,b_3,c_1,c_2,c_3).$$ 
%For such an algebra we associate the following matrix, $M_{{\mu},\{e_1, \cdots ,e_3\}}$,  belonging to $\mathcal{M}(3,9)$ the space of matrices of order $(3 \times 9)$ and given by
% $$\begin{pmatrix}
   %  \alpha_1v_1+\alpha_2v_2+\alpha_3v_3\\
    % \beta_1v_1+\beta_2v_2+\beta_3v_3 \\
     %\gamma_1v_1+\gamma_2v_2+\gamma_3v_3
%\end{pmatrix}
 %$$
% where
% $$\begin{array}{l}
  % v_1=(-\beta_3,\alpha_3,0,\beta_2,-\alpha_2,0,-\beta_1,\alpha_1,0) ,   \\
   % v_2=(-\gamma_3,0,\alpha_3,\gamma_2,0,-\alpha_2,-\gamma_1,0,\alpha_1),\\
   % v_3=(    0,-\gamma_3,\beta_3,0,\gamma_2,-\beta_2,0,-\gamma_1,\beta_1).
%\end{array}
%$$

\begin{theorem}
Any $3$-dimensional algebra is a Hom-Lie algebra.
\end{theorem}
\pf Consider a basis $\{e_1,e_2,e_3\}$ of the algebra $(A,\mu)$ and let $(\alpha_i,\beta_i,\gamma_i)$, $i=1,2,3$ be its structure constants defined previously. 
%Let $f$ be in $gl(3,\K)$ and let us consider the associated vector
%$$v_f=(a_1,a_2,a_3,b_1,b_2,b_3,c_1,c_2,c_3).$$ 
 The endomorphism  $f\in gl(3,\K) $
 satisfies the Hom-Jacobi condition if and only if its corresponding vector $v_f=\ ^t (a_1,a_2,a_3,b_1,b_2,b_3,c_1,c_2,c_3)$ is in the Kernel of the matrix 
 $$M_\mu=\begin{pmatrix}
23.1 & 23.2  & 23.3 & 31.1 & 31.2 & 31.3 & 12.1 & 12.2 & 12.3\\
  \end{pmatrix}
$$
where we use the notation $ij.k$ in place of $\mu(\mu(e_i,e_i),e_k)$
and
$$23.1=\begin{pmatrix}
    -\alpha_1 \beta_3 -\alpha_2\gamma_3  \\
   -\beta_1 \beta_3 -\beta_2\gamma_3  \\
    -\gamma_1 \beta_3 -\gamma_2\gamma_3  \\
\end{pmatrix} , \  23.2=\begin{pmatrix}
    \alpha_1 \alpha_3 -\alpha_3\gamma_3  \\
   \beta_1\alpha_3 -\beta_3\gamma_3  \\
    \gamma_1 \alpha_3 -\gamma_3\gamma_3  \\
\end{pmatrix} , \  23.3=\begin{pmatrix}
    \alpha_2 \alpha_3 +\alpha_3\beta_3  \\
   \beta_2\alpha_3 +\beta_3\beta_3  \\
    \gamma_2 \alpha_3 +\gamma_3\beta_3  \\
    \end{pmatrix} 
$$
$$31.1=\begin{pmatrix}
    \alpha_1 \beta_2 +\alpha_2\gamma_2  \\
   \beta_1 \beta_2 +\beta_2\gamma_2 \\
    \gamma_1 \beta_2 +\gamma_2\gamma_2  \\
\end{pmatrix} , \  31.2=\begin{pmatrix}
    -\alpha_1 \alpha_2 +\alpha_3\gamma_2  \\
   -\beta_1\alpha_2 +\beta_3\gamma_2  \\
    -\gamma_1 \alpha_2 +\gamma_3\gamma_2  \\
\end{pmatrix} , \  31.3=\begin{pmatrix}
    -\alpha_2 \alpha_2 -\alpha_3\beta_2  \\
   -\beta_2\alpha_2 -\beta_3\beta_2  \\
    -\gamma_2 \alpha_2 -\gamma_3\beta_2  \\
    \end{pmatrix} 
$$
$$12.1=\begin{pmatrix}
    -\alpha_1 \beta_1 -\alpha_2\gamma_1  \\
   -\beta_1 \beta_1 -\beta_2\gamma_1  \\
    -\gamma_1 \beta_1 -\gamma_2\gamma_1  \\
\end{pmatrix} , \  12.2=\begin{pmatrix}
    \alpha_1 \alpha_1 -\alpha_3\gamma_1  \\
   \beta_1\alpha_1 -\beta_3\gamma_1  \\
    \gamma_1 \alpha_1 -\gamma_3\gamma_1  \\
\end{pmatrix} , \  12.3=\begin{pmatrix}
    \alpha_2 \alpha_1 +\alpha_3\beta_1  \\
   \beta_2\alpha_1 +\beta_3\beta_1  \\
    \gamma_2 \alpha_1 +\gamma_3\beta_1  \\
    \end{pmatrix} .
$$Since the matrix $M_{\mu}$  belongs to $\mathcal{M}(3,9)$ and represents a linear morphism
 $$t: \K^9 \rightarrow \K^3.$$
 From the rank theorem we have 
 $$\dim {\rm Ker}\  t = 9- \dim {\rm Im}\  t \geq 6.$$
 Then this kernel is always non trivial and for any algebra $\mu$, there exists a non trivial element in the kernel. Then this algebra always admits a non trivial Hom-Lie structure.
 
 \medskip
 
 \noindent{\bf Consequence: On the classification of Hom-Lie algebras of dimension $3$} In \cite{ERdim3} one gives the classification of $3$-dimensional skew-symmetric $\K$-algebra. Since any Hom-Lie algebra is skew-symmetric, we deduce the classification of Hom-Lie algebras.

\medskip

We have seen that the kernel of $M_\mu$ is of dimension greater or equal to $3$. To compute this kernel we present in a vectorial form this matrix, using the notation $ij.k$ in place of $\mu(\mu(e_i,e_i),e_k)$. We have
$$M_\mu=\begin{pmatrix}
23.1 & 23.2  & 23.3 & 31.1 & 31.2 & 31.3 & 12.1 & 12.2 & 12.3\\
  \end{pmatrix}
$$
with
$$23.1=\begin{pmatrix}
    -\alpha_1 \beta_3 -\alpha_2\gamma_3  \\
   -\beta_1 \beta_3 -\beta_2\gamma_3  \\
    -\gamma_1 \beta_3 -\gamma_2\gamma_3  \\
\end{pmatrix} , \  23.2=\begin{pmatrix}
    \alpha_1 \alpha_3 -\alpha_3\gamma_3  \\
   \beta_1\alpha_3 -\beta_3\gamma_3  \\
    \gamma_1 \alpha_3 -\gamma_3\gamma_3  \\
\end{pmatrix} , \  23.3=\begin{pmatrix}
    \alpha_2 \alpha_3 +\alpha_3\beta_3  \\
   \beta_2\alpha_3 +\beta_3\beta_3  \\
    \gamma_2 \alpha_3 +\gamma_3\beta_3  \\
    \end{pmatrix} 
$$
$$31.1=\begin{pmatrix}
    \alpha_1 \beta_2 +\alpha_2\gamma_2  \\
   \beta_1 \beta_2 +\beta_2\gamma_2 \\
    \gamma_1 \beta_2 +\gamma_2\gamma_2  \\
\end{pmatrix} , \  31.2=\begin{pmatrix}
    -\alpha_1 \alpha_2 +\alpha_3\gamma_2  \\
   -\beta_1\alpha_2 +\beta_3\gamma_2  \\
    -\gamma_1 \alpha_2 +\gamma_3\gamma_2  \\
\end{pmatrix} , \  31.3=\begin{pmatrix}
    -\alpha_2 \alpha_2 -\alpha_3\beta_2  \\
   -\beta_2\alpha_2 -\beta_3\beta_2  \\
    -\gamma_2 \alpha_2 -\gamma_3\beta_2  \\
    \end{pmatrix} 
$$
$$12.1=\begin{pmatrix}
    -\alpha_1 \beta_1 -\alpha_2\gamma_1  \\
   -\beta_1 \beta_1 -\beta_2\gamma_1  \\
    -\gamma_1 \beta_1 -\gamma_2\gamma_1  \\
\end{pmatrix} , \  12.2=\begin{pmatrix}
    \alpha_1 \alpha_1 -\alpha_3\gamma_1  \\
   \beta_1\alpha_1 -\beta_3\gamma_1  \\
    \gamma_1 \alpha_1 -\gamma_3\gamma_1  \\
\end{pmatrix} , \  12.3=\begin{pmatrix}
    \alpha_2 \alpha_1 +\alpha_3\beta_1  \\
   \beta_2\alpha_1 +\beta_3\beta_1  \\
    \gamma_2 \alpha_1 +\gamma_3\beta_1  \\
    \end{pmatrix} .
$$
Thus, for a given algebra, it is easy to compute all the endomorphisms which satisfy the Jacobi-Hom-Lie condition.  Let us note also that the identity whose associated vector is 
$$v_{Id}=(1,0,0,0,1,0,0,0,1)$$ is in the kernel of $M_\mu$ if and only if $\mu$ satisfies
$$(23)1+(31)2+(12)3=0$$
that is if it is a Lie algebra. Let us note also, that the dimension of $M_\mu$ is an invariant associated with $\mu$, that is, if two Hom-Lie algebras are isomorphic then the Kernels of the associated matrices are isomorphic.

 \section{Dimension $4$}
  Let us resume the approach used for the study of the dimension $3$. 
  Let $(\mathcal{A},\mu)$ a $4$-dimensional $\K$-algebra.  Let us choose a basis $\{e_1,e_2,e_3,e_4\}$ of $\mathcal{A}$ and let us consider the corresponding constants structure of $\mu$:
  $$
  \left\{
  \begin{array}{l}
     \mu(e_1,e_2)= \alpha_1e_1+\beta_1e_2+\gamma_1e_3+\delta_1e_4   \\
     \mu(e_1,e_3	)= \alpha_2e_1+\beta_2e_2+\gamma_2e_3+\delta_2e_4   \\
     \mu(e_1,e_4)= \alpha_3e_1+\beta_3e_2+\gamma_3e_3+\delta_3e_4   \\
     \mu(e_2,e_3	)= \alpha_4e_1+\beta_4e_2+\gamma_4e_3+\delta_4e_4   \\
        \mu(e_2,e_4)= \alpha_5e_1+\beta_5e_2+\gamma_5e_3+\delta_5e_4   \\
     \mu(e_2,e_4	)= \alpha_6e_1+\beta_6e_2+\gamma_6e_3+\delta_6e_4   \\
\end{array}
\right.
$$
These algebra is a Hom-Lie algebra if there exists a linear endomorphism $f$ satisfying the Hom-Lie Jacobi equations. The endomorphism $f$ is represented in the basis $\{e_1,e_2,e_3,e_4\}$ by a square matrix of order $4$. The corresponding  $v_f$  belongs to $\K^{16}$:
  $$v_f=(a_1,a_2,a_3,a_4,b_1,b_2,b_3,b_4,c_1,c_2,c_3,c_4,d_1,d_2,d_3,d_4).$$
 These Hom-Lie Jacobi conditions  appear here in the form of a linear system on these coefficients, and this linear system which contains $16$ equations. Then $f$ satisfies the Hom-Lie conditions if and only if the vector $^tv_f$ is in the kernel of the associated matrix $M_{HL}(\mu)$ which is a square matrix of order $16$. Its kernel is not trivial if and only if 
  $$\det (M_{HL}(\mu))=0.$$
  We deduce
  
  \begin{proposition}
  The set $\mathcal{HL}_4$ of $4$-dimensional $\K$-Hom-Lie algebras is provided with a structure of affine algebraic variety embedded in $\K^{24}$.
  \end{proposition}
  
In fact,   $\det (M_{HL}(\mu))=0$ is a polynomial equation whose variables are the $24$ constants $\alpha_i,
\beta_i,\gamma_i,\delta_i$, $i=1,2,3,4.$  Let us note that this equation is homogeneous of degree $16$. 

For any $(\mathcal{A},\mu) \in \mathcal{HL}_4$  we can define the vector space $\ker M_{HL}(\mu)$. We thus define a singular vector bundle
$K(\mathcal{HL}_4)$ whose fiber over  $(\mathcal{A},\mu)$ is $\ker M_{HL}(\mu)$. This fiber corresponds to the set of Hom-Lie structure which can be defined on a given $4$-dimensional algebra.
 
 \noindent{\bf Remark.} In dimension $3$, $\mathcal{HL}_4$ is the affine variety $\mathcal{SA}lg_3$ which is isomorphic to the affine space $\K^9$.  We have seen that, in this case, the fibers of $K(\mathcal{HL}_3)$ are vector spaces of dimension greater or equal to $6$.

 \medskip
 We assume now that $\K$ is algebraically closed. In the previous presentation, it remains a problem concerning the existence or not of algebra $(\mathcal{A},\mu)$ for which we have $$\det (M_{HL}(\mu)) \neq 0.$$ We  meet then difficulties of calculation because a determinant of order $16$ is not easy to treat. Then we shall simplify this matricial approach. Since we assume that $\K$ is algebraically closed, the endomorphism $f$ can be reduced in an appropriated basis on the form
 $$\begin{pmatrix}
     a_1 & b_1 & 0 & 0    \\
      0& b_2  & c_2 & 0\\
      0 & 0 & c_3 & d_3 \\
      0 & 0 & 0 & d_4
\end{pmatrix}
$$ 
As a result, the matrix is reduced to a matrix $M'_{HL}(\mu)$ of type $16 \times 7$. This matrix is 

$$\left(
\begin{array}{cccccccccccccccc}
23 .1 & 31. 1 & 31.2 & 12. 2 & 12. 3  & 0 & 0\\
24. 1 & 41. 1 & 41. 2 & 0 & 0 &   12. 3   & 12. 4\\
34. 1 & 0 & 0 & 41. 2 & 41. 3 &  13. 3   & 13. 4\\
0 & 34. 1 & 34. 2  &    42. 2 & 42. 3 & 23. 3   & 23. 4\\
  \end{array}
  \right)
$$

with
{\tiny
$$
23.1=\left(
\begin{array}{l}
 -\beta_4 \alpha_1 - \gamma_4 \alpha_2 - \delta_4 \alpha_3\\
 -\beta_4 \beta_1 - \gamma_4 \beta_2 - \delta_4 \beta_3\\
 -\beta_4 \gamma_1 - \gamma_4 \gamma_2 - \delta_4 \gamma_3 \\
 -\beta_4 \delta_1 - \gamma_4 \delta_2 - \delta_4 \delta_3\\
\end{array}
\right), \ 
31.1=\left(
\begin{array}{l}
\beta_2 \alpha_1 + \gamma_2 \alpha_2 + \delta_2 \alpha_3\\
\beta_2 \beta_1 + \gamma_2 \beta_2 + \delta_2 \beta_3\\
 \beta_2 \gamma_1 + \gamma_2 \gamma_2 + \delta_2 \gamma_3\\
\beta_2 \delta_1 + \gamma_2 \delta_2 + \delta_2 \delta_3\\
\end{array}
\right), \  
31.2=\left(
\begin{array}{l}
 -\alpha_2 \alpha_1 + \gamma_2 \alpha_4 + \delta_2 \alpha_5\\
 -\alpha_2 \beta_1 + \gamma_2 \beta_4 + \delta_2 \beta_5\\
-\alpha_2 \gamma_1 + \gamma_2 \gamma_4 + \delta_2 \gamma_5\\
 -\alpha_2 \delta_1 + \gamma_2 \delta_4 + \delta_2 \delta_5\\
\end{array}
\right)
$$
}
{\tiny
$$
12.2=\left(
\begin{array}{l}
 \alpha_1 \alpha_1 - \gamma_1 \alpha_4 - \delta_1 \alpha_5\\
 \alpha_1 \beta_1 - \gamma_1 \beta_4 - \delta_1 \beta_5\\
 \alpha_1 \gamma_1 - \gamma_1 \gamma_4 - \delta_1 \gamma_5\\
 \alpha_1 \delta_1 - \gamma_1 \delta_4 - \delta_1 \delta_5\\
\end{array}
\right), \ 
12.3=\left(
\begin{array}{l}
\alpha_1 \alpha_2 + \beta_1 \alpha_4 - \delta_1 \alpha_6\\
 \alpha_1 \beta_2 + \beta_1 \beta_4 - \delta_1 \beta_6\\
\alpha_1 \gamma_2 + \beta_1 \gamma_4 - \delta_1 \gamma_6\\
 \alpha_1 \delta_2 + \beta_1 \delta_4 - \delta_1 \delta_6\\
\end{array}
\right), \
24.1=\left(
\begin{array}{l}
 -\beta_5 \alpha_1 - \gamma_5 \alpha_2 - \delta_5 \alpha_3 \\
 -\beta_5 \beta_1 - \gamma_5 \beta_2 - \delta_5 \beta_3 \\
-\beta_5 \gamma_1 - \gamma_5 \gamma_2 - \delta_5 \gamma_3\\
 -\beta_5 \delta_1 - \gamma_5 \delta_2 - \delta_5 \delta_3\\
\end{array}
\right)$$
}
{\tiny
$$
41.1= \left(
\begin{array}{l}
\beta_3 \alpha_1 + \gamma_3 \alpha_2 + \delta_3 \alpha_3\\
 \beta_3 \beta_1 + \gamma_3 \beta_2 + \delta_3 \beta_3\\
 \beta_3 \gamma_1 + \gamma_3 \gamma_2 + \delta_3 \gamma_3\\
 \beta_3 \delta_1 + \gamma_3 \delta_2 + \delta_3 \delta_3\\
\end{array}
\right), \ 
41.2=\left(
\begin{array}{l}
-\alpha_3 \alpha_1 + \gamma_3 \alpha_4 + \delta_3 \alpha_5\\
-\alpha_3 \beta_1 + \gamma_3 \beta_4 + \delta_3 \beta_5\\
-\alpha_3 \gamma_1 + \gamma_3 \gamma_4 + \delta_3 \gamma_5\\
 -\alpha_3 \delta_1 + \gamma_3 \delta_4 + \delta_3 \delta_5\\
\end{array}
\right), \ 
12.4=\left(
\begin{array}{l}
 \alpha_1 \alpha_3 + \beta_1 \alpha_5 + \gamma_1 \alpha_6\\
\alpha_1 \beta_3 + \beta_1 \beta_5 + \gamma_1 \beta_6\\
\alpha_1 \gamma_3 + \beta_1 \gamma_5 + \gamma_1 \gamma_6\\
\alpha_1 \delta_3 + \beta_1 \delta_5 + \gamma_1 \delta_6\\
\end{array}
\right) 
$$
}
{\tiny
$$
34.1=\left(
\begin{array}{l}
-\beta_6 \alpha_1 - \gamma_6 \alpha_2 - \delta_6 \alpha_3\\
 -\beta_6 \beta_1 - \gamma_6 \beta_2 - \delta_6 \beta_3\\
  -\beta_6 \gamma_1 - \gamma_6 \gamma_2 - \delta_6 \gamma_3\\
 -\beta_6 \delta_1 - \gamma_6 \delta_2 - \delta_6 \delta_3\\
\end{array}
\right), \ 
41.3=\left(
\begin{array}{l}
 -\alpha_3 \alpha_2 - \beta_3 \alpha_4 + \delta_3 \alpha_6\\
 -\alpha_3 \beta_2 - \beta_3 \beta_4 + \delta_3 \beta_6\\
 -\alpha_3 \gamma_2 - \beta_3 \gamma_4 + \delta_3 \gamma_6\\
 -\alpha_3 \delta_2 - \beta_3 \delta_4 + \delta_3 \delta_6\\
\end{array}
\right), \ 
13.3=\left(
\begin{array}{l}
\alpha_2 \alpha_2 + \beta_2 \alpha_4 - \delta_2 \alpha_6\\
 \alpha_2 \beta_2 + \beta_2 \beta_4 - \delta_2 \beta_6\\
 \alpha_2 \gamma_2 + \beta_2 \gamma_4 - \delta_2 \gamma_6\\
 \alpha_2 \delta_2 + \beta_2 \delta_4 - \delta_2 \delta_6\\
\end{array}
\right),
$$
}
{\tiny
$$
13.4=\left(
\begin{array}{l}
\alpha_2 \alpha_3 + \beta_2 \alpha_5 + \gamma_2 \alpha_6\\
 \alpha_2 \beta_3 + \beta_2 \beta_5 + \gamma_2 \beta_6\\
 \alpha_2 \gamma_3 + \beta_2 \gamma_5 + \gamma_2 \gamma_6\\
 \alpha_2 \delta_3 + \beta_2 \delta_5 + \gamma_2 \delta_6\\
\end{array}
\right), \  
34.2=
\left(
\begin{array}{l}
 \alpha_6 \alpha_1 - \gamma_6 \alpha_4 - \delta_6 \alpha_5\\
 \alpha_6 \beta_1 - \gamma_6 \beta_4 - \delta_6 \beta_5\\
 \alpha_6 \gamma_1 - \gamma_6 \gamma_4 - \delta_6 \gamma_5\\
 \alpha_6 \delta_1 - \gamma_6 \delta_4 - \delta_6 \delta_5\\
\end{array}
\right), \ 
42.2=
\left(
\begin{array}{l}
 -\alpha_5 \alpha_1 + \gamma_5 \alpha_4 + \delta_5 \alpha_5\\
 -\alpha_5 \beta_1 + \gamma_5 \beta_4 + \delta_5 \beta_5\\
-\alpha_5 \gamma_1 + \gamma_5 \gamma_4 + \delta_5 \gamma_5\\
 -\alpha_5 \delta_1 + \gamma_5 \delta_4 + \delta_5 \delta_5\\
\end{array}
\right) ,
$$
}
{\tiny
$$
42.3
\left(
\begin{array}{l}
-\alpha_5 \alpha_2 - \beta_5 \alpha_4 + \delta_5 \alpha_6\\
 -\alpha_5 \beta_2 - \beta_5 \beta_4 + \delta_5 \beta_6\\
 -\alpha_5 \gamma_2 - \beta_5 \gamma_4 + \delta_5 \gamma_6\\
 -\alpha_5 \delta_2 - \beta_5 \delta_4 + \delta_5 \delta_6\\
\end{array}
\right), \ 
23.3=\left(
\begin{array}{l}
 \alpha_4 \alpha_2 + \beta_4 \alpha_4 - \delta_4 \alpha_6\\
 \alpha_4 \beta_2 + \beta_4 \beta_4 - \delta_4 \beta_6\\
 \alpha_4 \gamma_2 + \beta_4 \gamma_4 - \delta_4 \gamma_6\\
 \alpha_4 \delta_2 + \beta_4 \delta_4 - \delta_4 \delta_6\\
\end{array}
\right), \  
23.4=\left(
\begin{array}{l}
 \alpha_4 \alpha_3 + \beta_4 \alpha_5 + \gamma_4 \alpha_6\\
 \alpha_4 \beta_3 + \beta_4 \beta_5 + \gamma_4 \beta_6\\
\alpha_4 \gamma_3 + \beta_4 \gamma_5 + \gamma_4 \gamma_6\\
 \alpha_4 \delta_3 + \beta_4 \delta_5 + \gamma_4 \delta_6\\
\end{array}
\right) $$
}
To a generic point $(\alpha_i,\beta_i,\gamma_i,\delta_i)$, $i=1,\cdots,6,$ this matrix is of rank $7$.  Then there exists a $4$-dimensional algebra
such that for any homomorphism $f$ the corresponding matrix written in a basis which reduces $f$ in a canonical form is of maximal rank. Such algebra cannot be a Hom-Lie algebra.
\begin{theorem}
The set  $\mathcal{HL}_4$ of the $4$-dimensional $\K$-Hom-Lie algebras is a affine algebraic variety strictly  contained in the affine plane $\mathcal{SA}lg_4$ isomorphic to $\K^{24}$.
  \end{theorem}
\pf From the previous calculus, $\mathcal{HL}_4$ is strictly contained in $\mathcal{SA}lg_4$. Since $\mathcal{HL}_4$ is determinate by the polynomial condition $\det (M_\mu)=0$,  we deduce that it is an affine algebraic sub variety of the affine plane $\mathcal{SA}lg_4$.

\medskip

Let us note that in the original basis $\{e_1,e_2,e_3,e_4\}$ the expression of $M_\mu$ is, using the notations of the previous section: 
{\tiny
\begin{equation}
\label{hl4}
\left(
\begin{array}{cccccccccccccccc}
23 .1 & 23.2  & 23. 3 & 23. 4 & 31. 1 & 31.2 & 31. 3 & 31. 4 &12. 1 & 12. 2 & 12. 3 & 12. 4 & 0 & 0 & 0 & 0\\
24. 1 & 24. 2  & 24. 3 & 24. 4 & 41. 1 & 41. 2 & 41. 3 & 41. 4 & 0 & 0 & 0 & 0 &12. 1 & 12. 2 & 12. 3   & 12. 4\\
34. 1 & 34. 2  & 34. 3 & 34. 4  & 0 & 0 & 0 & 0 & 41. 1 & 41. 2 & 41. 3 & 41. 4  &13. 1 & 13. 2 & 13. 3   & 13. 4\\
0 & 0 & 0 & 0 & 34. 1 & 34. 2  & 34. 3 & 34. 4  &  42. 1 & 42. 2 & 42. 3 & 42. 4  &23. 1 & 23. 2 & 23. 3   & 23. 4\\
  \end{array}
  \right)
\end{equation}
}
In a generic point the rank of this matrix is equal to $16$. Let us consider, for example, the following algebra:
$$
  \left\{
  \begin{array}{l}
     \mu(e_1,e_2)= e_2+2e_3-e_4   \\
     \mu(e_1,e_3	)= e_1+2e_2-e_3   \\
     \mu(e_1,e_4)= 2e_1-e_2+e_4   \\
     \mu(e_2,e_3	)= -e_1+e_3+2e_4   \\
     \mu(e_2,e_4	)= e_1+2e_2-e_3+3e_4   \\
     \mu(e_3,e_4)= -2e_1-e_2+e_3+2e_4   \\
\end{array}
\right.
$$
The associated matrix is 
$$
\left(
\begin{array}{cccccccccccccccc}
 -5 & -1 & 3 & -4 & -1 & 1 & 1 & -6 & 0 & 3 & -3 & -3 & 0 & 0 & 0 & 0 \\
 0 & -3 & 0 & 0 & 0 & -1 & -2 & -4 & 6 & 2 & -1 & 0 & 0 & 0 & 0 & 0 \\
 1 & 3 & -1 & 1 & 5 & -3 & -1 & 3 & 0 & -3 & 2 & 1 & 0 & 0 & 0 & 0 \\
 -2 & -9 & -4 & 1 & -2 & -1 & -4 & -5 & -2 & -1 & 4 & 7 & 0 & 0 & 0 & 0 \\
 -5 & -4 & 4 & 6 & 2 & 1 & -5 & -3 & 0 & 0 & 0 & 0 & 0 & 3 & -3 & -3 \\
 3 & 3 & 4 & 4 & -2 & 0 & -5 & 4 & 0 & 0 & 0 & 0 & 6 & 2 & -1 & 0 \\
 -5 & 6 & 1 & -3 & -2 & -5 & 4 & -1 & 0 & 0 & 0 & 0 & 0 & -3 & 2 & 1 \\
 -1 & -8 & -3 & 5 & 2 & 5 & 4 & 1 & 0 & 0 & 0 & 0 & -2 & -1 & 4 & 7 \\
 -5 & -1 & 3 & -7 & 0 & 0 & 0 & 0 & 2 & 1 & -5 & -3 & 1 & -1 & -1 & 6 \\
 1 & -6 & -2 & -1 & 0 & 0 & 0 & 0 & -2 & 0 & -5 & 4 & 0 & 1 & 2 & 4 \\
 3 & -3 & -1 & 2 & 0 & 0 & 0 & 0 & -2 & -5 & 4 & -1 & -5 & 3 & 1 & -3 \\
 -3 & -6 & -6 & -3 & 0 & 0 & 0 & 0 & 2 & 5 & 4 & 1 & 2 & 1 & 4 & 5 \\
 0 & 0 & 0 & 0 & -5 & -1 & 3 & -7 & 5 & 4 & -4 & -6 & -5 & -1 & 3 & -4 \\
 0 & 0 & 0 & 0 & 1 & -6 & -2 & -1 & -3 & 5 & -4 & -4 & 0 & -3 & 0 & 0 \\
 0 & 0 & 0 & 0 & 3 & -3 & -1 & 2 & 5 & -6 & -1 & 3 & 1 & 3 & -1 & 1 \\
 0 & 0 & 0 & 0 & -3 & -6 & -6 & -3 & 1 & 8 & 3 & -5 & -2 & -9 & -4 & 1 \\
\end{array}
\right)
$$
and its determinant is not zero. This algebra is not Hom-Lie.

\begin{corollary}
In the affine plane $\mathcal{SA}lg_4$, there is a Zariski open set whose elements are $4$-dimensional algebras without Hom-Lie structure.
\end{corollary}

\section{The general case}
We assume in tis section that $(\mathcal{A},\mu)$ is a $n$-dimensional $\K$-algebra.  Let $f$ be a linear endomorphism of $\mathcal{A}$. We fix a basis $\{e_1,\cdots,e_n\}$ of $\mathcal{A}$ and we consider the associated vector
$$v_f=(a_{1,1},\cdots,a_{n,1},a_{1,2},\cdots,a_{n,2},\cdots,a_{1,n},\cdots,a_{n,n}).$$
This vector satisfies the Jacobi-Hom-Lie condition if and only if it is in the kernel of the matrix $M_\mu$ where $M_\mu$ is the matrix $n^2$-columns and $\frac{n^2(n-1)(n-2)}{6}$ lines .  This matrix can be written 
{\tiny
$$\left(
\begin{array}{cccccccccccccccc}
23 .1 & \cdots & 23. n & 31. 1 &  \cdots & 31. n &12. 1 & \cdots & 12. n & 0 & 0 & 0 & 0 & \cdots  & \cdots  & 0   \\
24. 1 &  \cdots  & 24. n & 41. 1 & \cdots  & 41. 4 & 0 & \cdots & 0 &12. 1  & \cdots   & 12. 4 & 0 &\cdots   & \cdots & 0  \\
\cdots  &\cdots   & \cdots  & \cdots  & \cdots & \cdots  & \cdots  &\cdots  & \cdots  &\cdots  & \cdots   & \cdots & \cdots  & \cdots  & \cdots  & \cdots \\
34. 1  & \cdots  & 34. n  & 0 &  \cdots & 0 & 41. 1  & \cdots & 41. n  & 13. 2 & \cdots   & 13. n & 0 & \cdots  & \cdots  &0\\
0 &\cdots   & 0 & 34. 1  & \cdots & 34. n  &  42. 1 &\cdots  & 42. n  &23. 1 & \cdots   & 23. n & 0 & \cdots  & \cdots  & 0\\
\cdots  &\cdots   & \cdots  & \cdots  & \cdots & \cdots  & \cdots  &\cdots  & \cdots  &\cdots  & \cdots   & \cdots & \cdots  & \cdots  & \cdots  & \cdots \\
0  &\cdots   & 0  & 0  & \cdots & 0  &0  &\cdots  & 0 &0 & \cdots   & 0 & \cdots  & \cdots  & \cdots  & n-2 \ n-1.n \\

  \end{array}
  \right)
$$
}
 where $ij.k$ is the matrix $n \times 1$ whose line $l$ is
 $$\sum_s C_{i,j}^s C_{s,k}^l$$
 the constant $C_{i,j}^k$ being the structure constants of $\mu$ in the given basis $\{e_1,\cdots,e_n\}$.
 The kernel of this matrix is not trivial if and only if its rank is smaller or equal to $n^2-1$. This is equivalent to say that all the minor of order $n^2$ are zero.  We obtain a  system of polynomial equations on the variables $C_{i,j}^k$ of degree $n^4$. The number of equation corresponds to the number of combinations of $n^2$-elements amongst $\frac{n^2(n-1)(n-2)}{6}$ elements. For example, if $n=5$, the matrix $M_\mu$ is of order $25 \times 50$ and we have $\ds \frac{50!}{25!25!}$-polynomial equations of degree $625$ on the $50$ variables $C_{i,j}^k$.But this system is not reduced. For example, if we consider the first minor constituted with the first $25$ lines. If it is equal to zero, we can assume that the first line is a combination of the other $24$ lines. Then if we assume that this line is zero, all the minor containing this line are necessarily zero. This show that this system can be reduced to $\frac{n^2(n-1)(n-2)}{12}$ equations as soon as $n \geq 5$. 
 
 This system contains always a non trivial solutions, because Lie algebras are a special class of Lie algebra defined by considering that the identity is in the kernel of $M_\mu$.
 
 \begin{proposition}
 The set  $\mathcal{HL}_n$ of the $n$-dimensional $\K$-Hom-Lie algebras is a affine algebraic variety strictly  contained in the affine plane $\mathcal{A}lg_n$ isomorphic to $\K^{N}$ with $N=.\frac{n^2(n-1)(n-2)}{6}$
 \end{proposition}

\section{Some associated operads}

"An operad is an algebraic device, which encodes a type of algebras. Instead of studying the properties of a particular algebra, we focus on the universal operations that can be performed on the elements of any algebra of a given type. The information contained in an operad consists in these operations and all the ways of composing them" (Loday-Valette book)  It is illusory to think that there is an operad which encodes Hom-Lie algebras, because the latter include diverse sorts of algebras. However, as the operations occurring in the relation of Jacobi-Hom-Lie that is in the matrix $M_\mu$ are quadratic, we shall see  which kind of quadratic operas can be highlighted from the calculation of the pikernel of the matrix $M_\mu$. Let us remind briefly the definition of a quadratic  operad.

Let $\K[\sum_n]$ the $K$-algebrea of the symmetric group $\sum_n$. An operas $\mathcal{P}$ is defined by a sequence $\mathcal{P}(n), \ n \geq 1 $  of vector spaces over $\K$  such that $\mathcal{P}(n)$ is also a $\K[\sum_n]$-module and by composition maps: 
$$\circ_i :\mathcal{P}(n)\otimes \mathcal{P}(m)\rightarrow \mathcal{P}(n+m?1)$$
for  $i=1,...,n$ satisfying some associative properties called the May's axioms.

Any $\K[\sum_n]$-module $E$ generates a free operad denoted $\mathcal{F}(E)$ and satisfying 
$$\mathcal{F}(E)(1) = \K, \  \mathcal{F}(E)(2) = E. $$
In particular, if $E = \K[ \sum_2]$, the free module  $\mathcal{F} (E) (n)$ admits a basis whose elements are all the paranthezed products on $ n $ variables indexed by  par $\{1, 2, ..., n\} .$ 
Let $E$ be a  $\K[\sum_n]$-module and  $R$ a $\K[\sum_n]$-submodule of $\mathcal{F}(E)(3)$. On denote  by  $\mathcal{R}$ the ideal generated by $R$ that is the intersection of all the ideals $I$ of $\mathcal{F}(E)$ such that  $I(1) = {0}, I(2) = {0}$ and  $I(3) = R.$
We call binary quadratic operas generated by $E$ and the relations $R$, the operad $ \mathcal{P}(\K, E, R)$, also denoted by  $\mathcal{F}(E)/\mathcal{R}$  with
$$\mathcal{P}(\K, E, R)(n) = \mathcal{F}(E)(n)/\mathcal{R}(n).$$

Since we shall like determinate some quadratic operads based on Hom-Lie algebras, we assume that the $\sum_2$-module $E$ is generated by one element. and more precisely that $E ?= sgn_2$ where $ sgn_2$ is the one-dimensional signum representation. This is equivalent to say they the product is skew-symmetric. The ideal of relation $R$ will be given consdering special endomorphism $f$ in the kernel of $M_\mu$. For example, if we assume that $f = Id  \in \ker M_\mu$. This is implies that the coefficients of $M_\mu$ satisfy in dimension $3$ the following relation
$$
23.1 +31.2+12.3=0   \\
$$
and in dimension $4$, the relations
$$
\left\{
\begin{array}{l}
      23.1 +31.2+12.3=0   \\
      24.1+41.2 + 12.4 =0 \\
      34.1+41.3+13.4=0 \\
      34.2+42.3 +23.4 =0.
\end{array}
\right.
$$
This conduces to consider that $R$ is generated by the relation
$$ij.k + jk. i+ ki.j.$$ These relations are invariant by the action of $\sum_3$, this implies that $\mathcal{R}=R$ and the corresponding quadratic operad is the Lie operad $\mathcal{L}ie$.

\medskip

To find other quadratic operas in relayion with he Hom-Lie conditions, let us consider $\mathcal{F}(E)(3)$. This is a $12$-dimensional vector space generated by the triple
$$
\{ij.k, \ i.jk,  i\neq j \neq k, \ 1 \leq i,j,k \leq 3 \}.
$$
We consider that the multiplication is skew-symmetric, that is $\mathcal{R}(2)=\{ij+ji\}$. We deduce that $\mathcal{P}(3)=\mathcal{F}(E)(3)/\mathcal{R}(3)$ is a $4$-dimensional vector space generated by
$$
\{ij.k,  1 \leq i < j < k \leq 4 \}.
$$
Thus, the coefficients of the matrix ( \ref{hl4}) which belongs to $\mathcal{P}(3)$ are the elements of this last space. We deduce that a subclass of Hom-Lie algebras is associated with a quadratic operas if and only the equations
$$M_\mu (v_f)=0$$
constitue a linear system in $\mathcal{P}(3)$.  We have seen above that the vector $v_{Id}$ satisfies this condition. We have to determine all the others. From the form of the matrix ( \ref{hl4}), we deduce that 
$$v_f=(a_1,0,0,0,0,b_2,0,0,0,0,c_3,0,0,0,0,d_4)$$ and the equations of the kernel are
$$
\left\{
\begin{array}{l}
a_123 .1+b_2 31.2 +c_3 12. 3 =0\\
a_1 24. 1+b_2 41. 2+d_412. 4= 0\\
a_1 34. 1 +c_3 41. 3+d_4 13. 4=0\\
b_2 34. 2  +c_3 42. 3+d_4  23. 4=0\\
\end{array}
\right.
$$
This system is $\sum_3$-invariant if and only if
\begin{enumerate}
  \item $a_1=b_2=c_3=d_4=1$ and $f=Id$.
  \item  $ij.k=0.$
\end{enumerate}

The last case can be forgotten and we find only the $\mathcal{L}ie$ operad as binary quadratic operas in relation with the $Hom-Lie$ algebras.

\medskip

Since $M_\mu$ contains also termes of type $ij.i$, we are conduced to consider not only quadratic operas but also ternary operads. Then we consider $\mathcal{P}(2)=\mathcal{F}(E)(2)/\mathcal{R}(2)$ generated by $\{ij\}$, $\mathcal{P}(3)=\mathcal{F}(E)(3)/\mathcal{R}(3)$ generated by all the elements, that is $
\{ij.k,  1 \leq i < j < k \leq 4 \}.
$. The ideal of relations lies in $\mathcal{F}(E)(4)$. To illustrate this construction, let us take
$$v_f=(0,1,0,0,0,0,0,0,0)$$
Such ector is in the kernel of $M_\mu$ if and only if 
$$23.2=24.2=34.2=0.$$

\end{document}